\documentclass{article}
\usepackage{amsmath}
\usepackage{amsthm}
\usepackage{graphics}
\newtheorem{theorem}{Theorem}[section]
\newtheorem{lemma}[theorem]{Lemma}

\theoremstyle{definition}

\newtheorem{corollary}[theorem]{Corollary}
\newtheorem{proposition}[theorem]{Proposition}

\theoremstyle{remark}
\newtheorem{remark}[theorem]{Remark}
\newcommand{\eqb}{\begin{equation}}
\newcommand{\eqe}{\end{equation}}
\newcommand{\spb}{\begin{split}}
\newcommand{\spe}{\end{split}}
\newcommand{\cab}{\begin{cases}}
\newcommand{\cae}{\end{cases}}

\newcommand{\thmb}{\begin{theorem}}
\newcommand{\thme}{\end{theorem}}
\newcommand{\leb}{\begin{lemma}}
\newcommand{\lee}{\end{lemma}}
\newcommand{\cob}{\begin{corollary}}
\newcommand{\coe}{\end{corollary}}
\newcommand{\enb}{\begin{enumerate}}
\newcommand{\ene}{\end{enumerate}}
\newcommand{\mb}{\begin{matrix}}
\newcommand{\me}{\end{matrix}}
\newcommand{\cenb}{\begin{center}}
\newcommand{\cene}{\end{center}}

\newcommand{\prb}{\begin{proof}}
\newcommand{\pre}{\end{proof}}
\newcommand{\rb}{\begin{remark}}
\newcommand{\re}{\end{remark}}
\newcommand{\propb}{\begin{proposition}}
\newcommand{\prope}{\end{proposition}}
\newcommand{\lf}{\lfloor}
\newcommand{\rf}{\rfloor}
\newcommand{\lc}{\lceil}
\newcommand{\rc}{\rceil}
\newcommand{\lt}{\left}
\newcommand{\rt}{\right}
\newcommand{\la}{\langle}
\newcommand{\ra}{\rangle}
\newcommand{\R}{\mathbb{R}}

\newcommand{\Z}{\mathbb{Z}}
\newcommand{\I}{\mathbb{I}}
\newcommand{\N}{\mathbb{N}}
\newcommand{\cP}{\mathcal{P}}
\newcommand{\cT}{\mathbb{T}}
\newcommand{\cO}{\mathcal{O}}

\newcommand{\lala}{\langle\!\langle}
\newcommand{\rara}{\rangle\!\rangle}

\numberwithin{equation}{section}

\usepackage{amssymb}
\usepackage{amsfonts}
\usepackage{psfrag}
\usepackage{epsfig}
\usepackage{url}
\usepackage{mathtools}
\usepackage{extarrows}

\begin{document}

\title{Distributional asymptotics mod 1 of $(\log_bn)$}

 \author{Chuang Xu\\Department of Mathematical
   Sciences\\University of Copenhagen\\Copenhagen, {\sc
     Denmark}}
\maketitle

\begin{abstract}
\noindent
This paper studies the distributional asymptotics of the slowly changing sequence of logarithms $(\log_bn)$ with $b\in\N\setminus\{1\}.$ It is known that $(\log_bn)$ is not uniformly distributed modulo one, and its omega limit set is composed of a family of translated exponential distributions with constant $\log b.$ An improved upper estimate $\lt(\sqrt{\log N}/N\rt)$ is obtained for the rate of convergence with respect to (w.r.t.) the {\em Kantorovich metric} on the circle, compared to the general results on rates of convergence for a class of slowly changing sequences in the author's companion in-progress work. Moreover, a sharp rate of convergence $\lt(\log N/N\rt)$ w.r.t. the {\em Kantorovich metric} on the interval $[0,1]$, is derived. As a byproduct, the rate of convergence w.r.t. the {\em discrepancy metric} (or the {\em Kolmogorov metric}) turns out to be $\lt(\log N/N\rt)$ as well, which verifies that an upper bound for this rate derived in [{\sc Y. Ohkubo and O. Strauch}, {\em Distribution of leading digits of numbers}, Unif. Distrib. Theory, \textbf{11} (2016), no.1, 23--45.] is sharp.
\end{abstract}

\noindent
\hspace*{8.3mm}{\small {\bf Keywords.} Uniformly distributed modulo one sequence, slowly changing sequence, rate}\\
\hspace*{26.2mm} {\small  of convergence, Kantorovich metric, discrepancy, probability measures.}

\noindent
\hspace*{8.3mm}{\small {\bf MSC2010.} 11K06, 11K31, 60B10, 60E20.}

\section{Introduction}
Given a sequence of real numbers $(x_n),$ associate with it a sequence $(\nu_N(x_n))_{N\ge1}$ of finitely supported probability measures \[
\nu_N(x_n):=\frac{1}{N}\sum_{n=1}^N\delta_{\la x_n\ra},
\]
where $\delta_{\la x_n\ra}$ stands for the Dirac measure concentrated at  $\la x_n\ra$, the natural projection of $x_n$ onto the circle $\cT=\R/\Z$. Here and throughout, we write $\nu_N$ for $\nu_N(x_n)$ when $(x_n)$ is clear from the context. Note that $(\nu_N)$ is a sequence in the space $\cP(\cT)$ of all Borel probability measures on $\cT$. As a set, $\cP(\cT)$ can be identified with the subspace $\lt\{\mu\in\cP(\I): \mu(\{1\})=0\rt\}$ of $\cP(\I)$, where $\I$ denotes the compact unit interval $[0,1]$. Lowercase Greek letters $\mu$, $\nu$ are used henceforth to denote elements of both $\cP(\cT)$ and $\cP(\I)$, but it will always be clear from the context which space of measures is meant. Recall that a sequence $(x_n)$ in $\mathbb{R}$ is {\em uniformly distributed modulo one} (u.d.\ mod 1) \cite[Ch.1]{KN06} if
$\nu_N$ converges weakly in $\cP(\cT)$ to the uniform distribution $\lambda_{\cT}$ on $\cT$. Denote $\lambda_{\I}$ the uniform distribution on $\I$, and let $d_{\sf K}$ denote the {\em discrepancy} (or {\em Kolmogorov}) metric on $\cP(\I)$, i.e. \[d_{\sf K}(\mu,\nu)=\sup\nolimits_{x\in\I}\lt|\mu([0,x])-\nu([0,x])\rt|,\quad \forall\ \mu,\nu\in\cP(\I).\]  Recall from \cite[Cor.2.1.1]{KN06} that $(x_n)$ is u.d.\ mod 1 if and only if $\lim_{N\to\infty}d_{\sf K}(\nu_N\circ \iota^{-1},\lambda_{\I})=0$, where $\iota:\cT\to\I$ is the natural inclusion; see Section \ref{ch6-sec1} for details. It is well known \cite[Cor.2.1.2\&Thm.2.2.2]{KN06} that $d_{\sf K}(\nu_N\circ \iota^{-1},\lambda_{\I})\ge\frac{1}{2N}$ for every positive integer $N$; in fact, given any $(x_n)$ there exists a constant $c>0$ such that $d_{\sf K}(\nu_N\circ \iota^{-1},\lambda_{\I})>c\log N/N$ for infinitely many $N$.

There is a vast literature on the estimation of discrepancy, especially for u.d.\ mod 1 sequences. For instance, for the sequence $(an)$, where  $a\in\R$ is irrational with bounded partial quotients, \cite[Thm.2.3.4]{KN06} asserts that \eqb\label{ch6-eq1}
d_{\sf K}(\nu_N\circ\iota^{-1},\lambda_{\I})=\cO(\log N/N),
\eqe  and \eqref{ch6-eq1} also holds for Van der Corput sequence \cite[Thm.2.3.5]{KN06}.  However, much less research seems to have been undertaken on sequences that are not u.d.\ mod 1, for example, on {\em slowly changing sequences} \cite{K73}.

Given a sequence $(x_n)$ in $\R$, an improved notion with regard to the distributional asymptotics of $x_n$ is the {\em Omega limit set} $\Omega[x_n]$, defined as \[\Omega[x_n]=\bigl\{\mu\in\cP(\cT): \nu_{N_k}(x_n)\xlongrightarrow{k\to\infty}\mu\ \text{weakly for some sequence}\ (N_k)\ \text{in}\ \N\bigr\}.\] It is not hard to see that $\Omega[x_n]$ is non-empty, closed and connected \cite{W97}. For sequences $(x_n)$ that are {\em slowly changing} in the sense that \[\lim_{n\to\infty}n(x_{n+1}-x_n)=\xi\in\R,\] it was shown in \cite{K73} that $(x_n)$ is not u.d.\ mod 1; moreover, the elements of $\Omega[x_n]$, have been described in terms of asymptotic distribution functions. Similar results for slowly changing sequences in the literature include logarithms of natural numbers or prime numbers, iterated logarithms, and monotone functions of prime numbers \cite{GS08,K73,MS13,KN06,O11,OS16,SB06,W97,W35}. As far as the author knows, however, there were virtually no results, in the case of slowly changing sequences, on the rate(s) of convergence for subsequences of $\lt(\nu_N\rt)$ to $\Omega[x_n]$, not even for very basic sequences such as $\lt(\log_bn\rt)$ with $b\in\N\setminus\{1\}$, prior to \cite{X16}. Only recently did the author learn that \cite{OS16} establishes an upper bound $\lt(\log N/N\rt)$ for the latter, as well as their asymptotic distribution functions. Even there, however, the sharpness of the bound $\lt(\log N/N\rt)$ remains obscure. This article aims at resolving this obscurity. Specifically, for sequences $\lt(\log_bn\rt)$, every limit point in $\Omega[x_n]$ is clearly identified, and $\lt(\log N/N\rt)$ is shown to be the sharp rate of convergence w.r.t. $d_{\sf K}$.

While the discrepancy metric (on $\cP(\cT)$, as induced by $d_{\sf K}$) has been used in uniform distribution theory for decades, its usage for sequences that are {\em not} u.d.\ mod 1 appears debatable. In fact, when analyzing such sequences, it may be more natural to study $\Omega[x_n]$ with a metric metrizing the weak topology of $\cP(\cT)$ such as, for instance, the {\em Kantorovich} (or {\em transport}) metric $d_{\cT}$. In a recent note \cite{X18}, the author obtained several results in this regard, including an upper bound  $\lt(\log N/N\rt)$ for $d_{\cT}$-convergence. As is shown in this article, however, this bound is not sharp, and better bounds are provided to replace it for $(x_n)=\lt(\log_bn\rt)$. From the arguments presented, it will also become evident that finding a good {\em lower} bound remains a formidable challenge, even for sequences as simple as $(\log_bn)$.

\section{Preliminaries and notations}\label{ch6-sec1}

Let $\R$, $\Z$, and $\N$ be the set of real numbers, integers, and positive integers, respectively. Recall that $\cT=\R/\Z$ can be thought of geometrically as the unit circle $\lt\{e^{2\pi ix}: x\in\R\rt\}$ in the complex plane, with its usual topology. For $-\infty<a<b\le\infty,$ let $[a,b[:=\{y\in\mathbb{R}:a\le y<b\}$; intervals $[a,b],\ ]a,b],\ ]a,b[$ are defined analogously. Let $\lf x\rf$, $\lc x\rc$, and $\lala x\rara=x-\lf x\rf$ be the floor (i.e., the largest integer $\le x$), the ceiling (i.e., the smallest integer $\ge x$), and the fractional part of $x\in\R$, respectively. In what follows, it will prove useful to denote by $\pi: \R\to\cT$, with $\pi(x)=\la x\ra=x+\Z$, and by $\iota: \cT\to\I$, with $\iota(\la x\ra)=\lala x\rara$, the natural projection and inclusion, respectively. Arguably the most fundamental maps on $\cT$ are rotations: Given any $\theta\in\R$, let $R_{\theta}$ be the (counter-clockwise) rotation of $\cT$ by $2\pi\theta$, that is, $R_{\theta}(\la x\ra)=\la x+\theta\ra$ for all $\la x\ra\in\cT$. With this, clearly $R^k_{\theta}=R_{k\theta}=R_{\lala k\theta\rara}$ for all $\theta\in\R$ and $k\in\Z$.

Let $(X,\rho_X)$ be a compact metric space, and $\cP(X)$ the space of all Borel probability measures on $X,$  endowed with the weak topology. Recall that $\cP(X)$ is compact and metrizable. The \textit{Kantorovich distance} on $X$ is
\[
d_X(\mu,\nu)=\inf\nolimits_{\gamma}\int_{X\times X}\rho_X(x,y){\rm d}\gamma(x,y),\quad \forall\ \mu,\nu\in\cP(X),
\]
where the infimum is taken over all Borel probability measures $\gamma$ on $X\times X$ with marginals $\mu$ and $\nu.$ Note that $d_X$ metrizes the weak topology of $\mathcal{P}_X$. For $X=\I$ and $X=\cT$, let $\rho_\I=|x-y|$, and $\rho_{\cT}(x,y)=\min\{|\iota(x)-\iota(y)|,1-|\iota(x)-\iota(y)|\}$, $\forall\ x,y\in X$, respectively. Note also that $\mu\mapsto\mu\circ\pi^{-1}$ maps $\cP(\I)$ continuously onto $\cP(\cT)$; when restricted to $\widetilde{\cP}(\I):=\lt\{\mu\in\cP(\I): \mu(\{1\})=0\rt\}$, a dense $G_{\delta}$-set in $\cP(\I)$, this even yields a continuous bijection, but not a homeomorphism, as $\widetilde{\cP}(\I)$ is not compact. In the opposite direction, $\mu\mapsto\mu\circ\iota^{-1}$ establishes a measurable bijection from $\cP(\cT)$ onto $\widetilde{\cP}(\I)$. Note also that $\mu\mapsto\mu\circ R_{\theta}^{-1}$ defines a homeomorphism of $\cP(\cT)$.

Recall that $d_{\I}$ can be expressed explicitly as
\eqb\label{ch6-eq2}
d_{\I}(\mu,\nu)=\int_0^1\lt|F_{\mu}(x)-F_{\nu}(x)\rt|{\rm d}x,\quad\, \forall\ \mu,\nu\in\cP(\I),
\eqe
where $F_{\nu}$ is the distribution function of $\nu$.
A method of computing $d_{\cT}$ has been developed in \cite{CM95}; only the following simple upper bound will be used here.

\propb\label{ch6-prop1}{\rm\cite[Cor.3.8]{CM95}.} Assume that $\mu,\ \nu\in\cP(\cT)$. Then
\[\begin{split} d_{\cT}(\mu,\nu)&\le\inf\nolimits_{y\in\I}\int_0^1\bigl|\lt(F_{\mu\circ\iota^{-1}
}(x)-F_{\nu\circ\iota^{-1}}(x)\rt)-\lt(F_{\mu\circ\iota^{-1}}(y)-F_{\nu\circ
\iota^{-1}}(y)\rt)\bigr|{\rm d}x\\
&\le d_{\sf K}\lt(\mu\circ\iota^{-1},\nu\circ\iota^{-1}\rt).\end{split}\]
\prope

For every $a>0$, consider the negative exponential distribution $-{\rm Exp}(a)$ on $\R$ with parameter $a$, that is,
\[F_{{\rm Exp}(a)}(x)=e^{ax},\quad\, \forall\ x\le0,\] and let $E_a=-{\rm Exp}(a)\circ\pi^{-1}\in\cP(\cT)$. Thus
\[F_{E_a\circ \iota^{-1}}(x)=\frac{e^{ax}-1}{e^a-1},\quad\, \forall\ x\in\I.\] Rotated versions of $E_a$, that is, probabilities $E_a\circ R_{\theta}^{-1}$ with $\theta\in\R$, play an important role in this article. For such probabilities, observe that \[F_{E_a\circ R_{\theta}^{-1}}\circ \iota^{-1}=\cab \frac{e^{\lala\theta\rara\lt(e^{ax}-1\rt)}}{e^a-1}\hfill\text{if}\ x\in[0,1-\lala\theta\rara],\\ 1+\frac{e^{\lala\theta\rara\lt(e^{a(x-1)}-1\rt)}}{e^a-1}\hfill\text{if}\ x\in[1-\lala\theta\rara,1[.\cae\]Henceforth, our analysis focuses on the sequences $(x_n)=\lt(\log_bn\rt)$ with $b\in\N\setminus\{1\}$, and the associated discrete measures $\nu_N=\nu_N\lt(\log_bn\rt)\in\cP(\cT)$. A simple calculation yields an explicit formula for the distribution function of $\nu_N\circ\iota^{-1}$.
\propb\label{ch6-prop2}
Assume that $b\in\N\setminus\{1\}$ and $N\in\N$. Then, with $L=\lf\log_bN\rf$,
\eqb\label{ch6-eq3}
\begin{split} &F_{\nu_N\circ\iota^{-1}}(x)\\
=&\cab \displaystyle\frac{L+1+\sum_{j=0}^L\lt(\lf ib^{-j}\rf-b^{L-j}\rt)}{N}\ \ \quad\quad\quad\quad\ \text{if}\ x\in\lt[\log_b\frac{i}{b^L},\log_b\frac{i+1}{b^L}\rt[,\ i=b^L,\ldots,N-1,\\
 1+\displaystyle\frac{L+1+\sum_{j=0}^L \lt(\lf\lf Nb^{-1}\rf b^{-j}\rf-b^{L-j}\rt)}{N}\ \text{if}\ x\in\biggl[\log_b\frac{N}{b^L},\log_b\frac{b\lf Nb^{-1}\rf+b}{b^L}\biggr[,\\ 1+\displaystyle\frac{L+1+\sum_{j=0}^L \lt(\lf ib^{-j}\rf-b^{L-j}\rt)}{N}\quad\quad\ \ \ \ \text{if}\ x\in\lt[\log_b\frac{bi}{b^L},\log_b\frac{b(i+1)}{b^L}\rt[,\ i=\lf Nb^{-1}\rf+1,\\ \ \quad\quad\quad\quad\quad\quad\quad\quad\quad\quad\quad\quad\quad\quad\quad\quad\quad\quad\quad\quad\quad\quad\quad\quad\quad\quad\quad\quad\quad
 \quad\quad\quad\quad\quad  \ldots,b^L-1.\cae\end{split}\eqe
\prope

\section{Rates of convergence}\label{ch6-sec2}
In this section, we study the rate of convergence for subsequences of $\lt(\nu_N\rt)_{N\ge1}$ w.r.t.\  $d_{\cT}$, $d_{\I}$, and $d_{\sf K}$. Throughout, for ease of exposition, all proofs are given for $b=10$, but all arguments can easily be adjusted to any other base $b\in\N\setminus\{1\}$.

\subsection{Upper bound for the rate of convergence w.r.t.\ $d_{\cT}$}\label{ch6-sec2-sub1}

We first present our main result regarding an upper bound for the rate of convergence w.r.t.\ $d_{\cT}$.
\thmb\label{ch6-thm1}
Assume $b\in\N\setminus\{1\}.$ Then
\eqb\label{ch6-eq4}
\limsup_{N\to\infty}\frac{N}{\sqrt{\log N}}d_{\mathbb{T}}\lt(\nu_N,E_{\log b}\circ R^{-1}_{-\log_bN}\rt)<+\infty.
\eqe
\thme
\prb
Recall that $b=10$ is assumed throughout; for every $N\in\N$, let $n=\lf\log_{10}N\rf+1$ for convenience, and thus $10^{n-1}\le N\le10^n-1$, as well as $\eta_N=E_{\log_bN}\in\cP(\cT)\circ R^{-1}_{-\log_bN}\in\cP(\cT)$. By Proposition \ref{ch6-prop1}, it suffices to estimate $$\int_0^1\Bigl|\bigl(F_{\nu_N\circ\iota^{-1}}(x)-F_{\eta_N\circ\iota^{-1}}(x)
\bigr)-\bigl(F_{\nu_N\circ\iota^{-1}}(y)
-F_{\eta_N\circ\iota^{-1}}(y)\bigr)\Bigr|{\rm d}x,$$ for an appropriate $0\le y<1$. Utilizing Proposition \ref{ch6-prop2} we first simplify the latter expression as follows: For every $y\in\lt[0,\log_{10}N-n+1\rt[,$ let $i_0=\lf10^{y+n-1}\rf.$ Then $10^{n-1}\le i_0\le N-1$ and $y\in\lt[\log_{10}i_0-n+1,\log_{10}(i_0+1)-n+1\rt[.$ Similarly, for $10^{n-1}\le i\le 10^n-1$ and $x\in\I$, \eqb\label{ch6-eq5}
x\in\lt[\log_{10}i-n+1,\log_{10}(i+1)-n+1\rt[\quad\Leftrightarrow\ \ \  i=\lf10^{x+n-1}\rf.\eqe With this, it follows from Proposition \ref{ch6-prop2} that
\begin{align*}
&\int_0^1\Bigl|\bigl(F_{\nu_N\circ\iota^{-1}}(x)-F_{\eta_N\circ\iota^{-1}}(x)
\bigr)-\bigl(F_{\nu_N\circ\iota^{-1}}(y)
-F_{\eta_N\circ\iota^{-1}}(y)\bigr)\Bigr|{\rm d}x\\
     =&\sum_{i=10^{n-1}}^{N-1}\int_{\log_{10}i-n+1}^{\log_{10}(i+1)-n+1}\lt|
     \frac{\sum_{j=0}^{n-1}\lt(\lf i10^{-j}\rf-\lf i_010^{-j}\rf\rt)}{N}-\frac{10^n\lt(10^x-10^y\rt)}{9N}\rt|{\rm d}x\\
     &+\int_{\log_{10}(N/10)-n+2}^{\log_{10}\lt(\lf N/10\rf+1\rt)-n+2}\biggl|\frac{\sum_{j=0}^{n-1}\lt(\lf \lf N/10\rf10^{-j}\rf-\lf i_010^{-j}\rf\rt)}{N}-\frac{10^n(10^{x-1}-10^y)}{9N}\biggr|{\rm d}x\\&+
     \sum_{i=\lf N/10\rf+1}^{10^{n-1}-1}\int_{\log_{10}i-n+2}^{\log_{10}(i+1)-n+2}\lt
     |\frac{\sum_{j=0}^{n-1}\lt(\lf i10^{-j}\rf-\lf i_010^{-j}\rf\rt)}{N}-\frac{10^n(10^{x-1}-10^y)}{9N}\rt|{\rm d}x\\
     =&\sum_{i=\lf N/10\rf+1}^{N-1}\int_{\log_{10}i-n+1}^{\log_{10}(i+1)-n+1}\lt|\frac{
     \sum_{j=0}^{n-1}\lt(\lf i10^{-j}\rf-\lf i_010^{-j}\rf\rt)}{N}-\frac{10^n(10^x-10^y)}{9N}\rt|{\rm d}x\\&+\int_{\log_{10}(N/10)-n+1}^{\log_{10}\lt(\lf N/10\rf+1\rt)-n+1}\biggl|\frac{\sum_{j=0}^{n-1}\lt(\lf \lf N/10\rf10^{-j}\rf-\lf i_010^{-j}\rf\rt)}{N}-\frac{10^n(10^x-10^y)}{9N}\biggr|{\rm d}x.
    \end{align*}
Since $\lt|F_{\nu_N\circ\iota^{-1}}(x)-F_{\eta_N\circ\iota^{-1}}(x)\rt|\le1$ for all $x\in\I$, it easily follows that
\[
  \int_{\log_{10}(N/10)-n+1}^{\log_{10}(\lf N/10\rf+1)-n+1}\biggl|\frac{\sum_{j=0}^{n-1}(\lf \lf N/10\rf10^{-j}\rf-\lf i_010^{-j}\rf)}{N}-\frac{10^n(10^x-10^y)}{9N}\biggr|{\rm d}x=\cO\lt(N^{-1}\rt).
\]From \eqref{ch6-eq5} and $i_0=\lf10^{y+n-1}\rf$, it is readily verified that
\begin{align*}
  &\frac{\sum_{j=0}^{n-1}(\lf i10^{-j}\rf-\lf i_010^{-j}\rf)}{N}-\frac{10^n(10^x-10^y)}{9N} \\ =&\frac{\sum_{j=0}^{n-1}\lt((i_010^{-j}-\lf i_010^{-j}\rf)-(i10^{-j}-\lf i10^{-j}\rf)\rt)}{N}-\frac{10^{1-n}(\lf10^{x+n-1}\rf-\lf10^{y+n-1}
  \rf)}{9N}\\
  &+\frac{10\lt((\lf10^{x+n-1}\rf-10^{x+n-1})-(\lf10^{y+n-1}\rf-10
  ^{y+n-1})\rt)}{9N}.
\end{align*}
Since also \[\begin{split}
&  \lt|\frac{10\lt((\lf10^{x+n-1}\rf-10^{x+n-1})-(\lf10^{y+n-1}\rf-
10^{y+n-1})\rt)}{9N}-\frac{10^{1-n}(\lf10^{x+n-1}\rf-\lf10^{y+n-1}
\rf)}{9N}\rt|\\=&\cO\lt(N^{-1}\rt),\quad \text{for}\ y\in\lt[0,\log_{10}N-n+1\rt[,
\end{split}
\]
 we obtain \eqb\label{ch6-eq6}
\begin{split}&\int_0^1\Bigl|\bigl(F_{\nu_N\circ\iota^{-1}}(x)-F_{\eta_N\circ\iota^{-1}}(x)
\bigr)-\bigl(F_{\nu_N\circ\iota^{-1}}(y)
-F_{\eta_N\circ\iota^{-1}}(y)\bigr)\Bigr|{\rm d}x\\
     =&\frac{1}{N}\sum_{i=\lf N/10\rf+1}^{N-1}\log_{10}\lt(1+\frac{1}{i}\rt)\cdot\lt|\sum_{j=0}^{n-1}
     \lt((i_010^{-j}-\lf i_010^{-j}\rf)-(i10^{-j}-\lf i10^{-j}\rf)\rt)\rt|+\cO\lt(N^{-1}\rt).
     \end{split}
    \eqe
In the following, we further estimate the right hand side of \eqref{ch6-eq6}. The elementary inequality $$x-x^2/2\le\log(x+1)\le x,\quad \forall\ x\ge0$$ yields $$\frac{1}{i\log10}-\frac{1}{2i^2\log10}\le\log_{10}\lt(1+\frac{1}{i}\rt)\le
\frac{1}{i\log10},$$
and we also have\[\begin{split}&\frac{1}{N}\sum_{i=\lf N/10\rf+1}^{N-1}\frac{1}{2i^2\log10}\lt|\sum_{j=0}^{n-1}\lt((i_010^{-j}-\lf i_010^{-j}\rf)-(i10^{-j}-\lf i10^{-j}\rf)\rt)\rt|=\cO\lt(N^{-2}\log N\rt). \end{split}\]Hence
\begin{align*}
     &\int_0^1\Bigl|\bigl(F_{\nu_N\circ\iota^{-1}}(x)-F_{\eta_N\circ\iota^{-1}}(x)
\bigr)-\bigl(F_{\nu_N\circ\iota^{-1}}(y)
-F_{\eta_N\circ\iota^{-1}}(y)\bigr)\Bigr|{\rm d}x\\
     =&\frac{1}{N\log10}\sum_{i=\lf N/10\rf+1}^{N-1}\frac{1}{i}\lt|\sum_{j=0}^{n-1}\lt(i_010^{-j}-\lf i_010^{-j}\rf)-(i10^{-j}-\lf i10^{-j}\rf)\rt)\rt|+\cO\lt(N^{-1}\rt).
    \end{align*}
     From\[\frac{1}{N\log10}\frac{1}{i}\lt|\sum_{j=0}^{n-1}\lt((i_010^{-j}-\lf i_010^{-j}\rf)-(i10^{-j}-\lf i10^{-j}\rf)\rt)\rt|\le\frac{2n}{N\lf N/10\rf\log10},\] for $i=\lf N/10\rf$ and $i=N$, it follows that \eqb\label{ch6-eq7}
     \begin{split}
     &\int_0^1\Bigl|\bigl(F_{\nu_N\circ\iota^{-1}}(x)-F_{\eta_N\circ\iota^{-1}}(x)
\bigr)-\bigl(F_{\nu_N\circ\iota^{-1}}(y)
-F_{\eta_N\circ\iota^{-1}}(y)\bigr)\Bigr|{\rm d}x\\
     =&\frac{1}{N\log10}\sum_{i=\lf N/10\rf}^N\frac{1}{i}\lt|\sum_{j=0}^{n-1}\lt((i_010^{-j}-\lf i_010^{-j}\rf)-(i10^{-j}-\lf i10^{-j}\rf)\rt)\rt|+\cO\lt(N^{-1}\rt).
    \end{split}\eqe
         Completely analogous arguments show that \eqref{ch6-eq7} holds also for $y\in[\log_{10}N-n+1,1[$ with $i_0=\lf10^{y+n-2}\rf.$ Thus it suffices to determine the constant order of amplitude of
$$\sum_{i=\lf N/10\rf}^N\frac{1}{i}\lt|\sum_{j=0}^{n-1}\lt((i_010^{-j}-\lf i_010^{-j}\rf)-(i10^{-j}-\lf i10^{-j}\rf)\rt)\rt|.$$
To get rid of the absolute value, one can use the Cauchy-Schwarz inequality: \eqb\label{ch6-eq8}\begin{split}
     &\lt\{\sum_{i=\lf N/10\rf}^N\frac{1}{i}\lt|\sum_{j=0}^{n-1}\lt((i_010^{-j}-\lf i_010^{-j}\rf)-(i10^{-j}-\lf i10^{-j}\rf)\rt)\rt|\rt\}^2\\
     \le&\sum_{i=\lf N/10\rf}^N\frac{1}{i^2}\sum_{i=\lf N/10\rf}^N\lt\{\sum_{j=0}^{n-1}\lt((i_010^{-j}-\lf i_010^{-j}\rf)-(i10^{-j}-\lf i10^{-j}\rf)\rt)\rt\}^2.
    \end{split}\eqe Note that
    \eqb\label{ch6-eq9}
    \sum_{i=\lf N/10\rf}^N\frac{1}{i^2}=9N^{-1}+\cO(N^{-2}).\eqe It remains to estimate $$\sum_{i=0}^N\Bigl(\sum_{j=0}^{n-1}\bigl((i_010^{-j}-\lf i_010^{-j}\rf)-(i10^{-j}-\lf i10^{-j}\rf)\bigr)\Bigr)^2,$$ which can be rewritten as \eqb\label{ch6-eq10}
    \begin{split}
     &\sum_{i=0}^N\lt\{\sum_{j=0}^{n-1}\lt((i_010^{-j}-\lf i_010^{-j}\rf)-(i10^{-j}-\lf i10^{-j}\rf)\rt)\rt\}^2\\=&(N+1)\lt(\sum_{j=0}^{n-1}(i_010^{-j}-\lf i_010^{-j}\rf)\rt)^2
     +\sum_{i=0}^N\lt(\sum_{j=0}^{n-1}(i10^{-j}-\lf i10^{-j}\rf)\rt)^2\\&-2\sum_{j=0}^{n-1}(i_010^{-j}-\lf i_010^{-j}\rf)\sum_{i=0}^N\sum_{j=0}^{n-1}(i10^{-j}-\lf i10^{-j}\rf).
    \end{split}\eqe
   In the following, we consider each term on the right-hand side of \eqref{ch6-eq10} individually.

First we consider $\sum_{i=0}^N\sum_{j=0}^{n-1}(i10^{-j}-i\lf 10^{-j}\rf),$
by switching the order of the summations.
    For every $i=0,\cdots,N$ and $j=0,\cdots,n-1,$ there exist nonnegative integers $k,\ l$ with $l\le10^j-1$ such that $i=k10^j+l,$ and hence $i10^{-j}-\lf i10^{-j}\rf=l10^{-j}.$ Therefore\begin{align*}\lt\{i:0\le i\le N\rt\}=&\lt\{k10^j+l: 0\le k\le\lf N10^{-j}\rf-1,  0\le l\le10^j-1\rt\}\\&\cup\lt\{\lf N10^{-j}\rf10^j+l: 0\le l\le N-\lf N10^{-j}\rf10^j\rt\}.\end{align*} Let  $N=\overline{a_{n-1}\cdots a_0}=\sum_{j=0}^{n-1}a_j10^j$ with $0\le a_j\le 9$ for all $0\le j\le n-1$.
Notice that \eqb\label{ch6-eq11}
\lf N10^{-j}\rf=N10^{-j}-\sum_{r=0}^{j-1}a_r10^{r-j},\quad \forall\ j=1,\cdots,n-1.\eqe From the simple observations $$\sum_{j=0}^{n-1}(1-10^{-j})=n+\cO(1)\,\quad \text{and}\,\quad \sum_{j=1}^{n-1}10^{-j}\sum_{r=0}^{j-1}a_r10^r=\cO(n),$$ it is tedious but straightforward to deduce that
    \eqb\label{ch6-eq12}
   \sum_{i=0}^N\sum_{j=0}^{n-1}\lt(i10^{-j}-\lf i10^{-j}\rf\rt)
      =\frac{1}{2}\lt(n-\frac{10}{9}\rt)N-\frac{1}{2}\sum_{j=1}^{n-1}
      \sum_{r=0}^{j-1}a_r10^r+\frac{1}{2}\sum_{j=1}^{n-1}10^{-j}\lt(
      \sum_{r=0}^{j-1}a_r10^r\rt)^2+\cO(n).\eqe
Next, we deal with  $\sum_{i=0}^N\lt(\sum_{j=0}^{n-1}(i10^{-j}-\lf i10^{-j}\rf)\rt)^2,$ which can be expanded as
    \begin{align}\label{ch6-eq13}
    &\sum_{i=0}^N\lt(\sum_{j=0}^{n-1}(i10^{-j}-\lf i10^{-j}\rf)\rt)^2\notag\\
    =&2\sum_{i=0}^N\sum_{j=1}^{n-1}(i10^{-j}-\lf i10^{-j}\rf)\sum_{r=0}^{j-1}(i10^{-r}-\lf i10^{-r}\rf)+\sum_{i=0}^N\sum_{j=0}^{n-1}(i10^{-j}-\lf i10^{-j}\rf)^2.
    \end{align}
For every $1\le j\le n-1$, let $K_j=N-\lf N10^{-j}\rf10^j$ for notational convenience. Then similarly, \begin{align*}
  \lt\{i:0\le i\le N\rt\}=&\lt\{k10^j+p10^r+l: 0\le k\le\lf N10^{-j}\rf-1, 0\le p\le10^{j-r}-1, 0\le l\le10^r-1\rt\}\\&\cup\lt\{\lf N10^{-j}\rf10^j+p10^r+l: 0\le p\le\lf K_j10^{-r}\rf-1, 0\le l\le10^r-1
    \rt\}\\&\cup\lt\{\lf N10^{-j}\rf10^j+\lf K_j10^{-r}\rf10^r+l: 0\le l\le K_j-\lf K_j10^{-r}\rf10^r\rt\},
\end{align*}
from which it follows that\eqb\label{ch6-eq14}
\begin{split}
  &\sum_{i=0}^N\sum_{j=1}^{n-1}(i10^{-j}-\lf i10^{-j}\rf)\sum_{r=0}^{j-1}(i10^{-r}-\lf i10^{-r}\rf)\\
  =&\sum_{j=1}^{n-1}\sum_{r=0}^{j-1}\lt(\sum_{k=0}^{\lf N10^{-j}\rf-1}\sum_{p=0}^{10^{j-r}-1}\sum_{l=0}^{10^r-1}(p10^r+l)\cdot10^
  {-j}l10^{-r}+\sum_{p=0}^{\lf K_j10^{-r}\rf-1}\sum_{l=0}^{10^r-1}(p10^r+l)10^{-j}l10^{-r}\rt.\\&\lt.\quad\quad\quad\quad+\sum_{l=0}^
  {K_j-\lf K_j10^{-r}\rf10^r}\lt(\lf K_j10^{-r}\rf10^r+l\rt)10^{-j}l10^{-r}\rt).
\end{split}
\eqe
From \eqref{ch6-eq11} and \eqref{ch6-eq14}, a lengthy but elementary calculation leads to \begin{align*}
  &\sum_{i=0}^N\sum_{j=0}^{n-1}(i10^{-j}-\lf i10^{-j}\rf)\sum_{r=0}^{j-1}(i10^{-r}-\lf i10^{-r}\rf)\\
  =&\lt(\frac{n^2}{8}-\frac{85n}{216}\rt)N-\frac{1}{4}\sum_{j=1}^{n-1}j
  \sum_{l=0}^{j-1}a_l10^l+\frac{1}{4}\sum_{j=1}^{n-1}j10^{-j}\lt(\sum_{l=0}^
  {j-1}a_l10^l\rt)^2+\cO(N).\end{align*}
Analogously, one obtains also
$\sum_{i=0}^N\sum_{j=0}^{n-1}(i10^{-j}-\lf i10^{-j}\rf)^2=\frac{nN}{3}+\cO(N).$
Note that \eqref{ch6-eq13} immediately leads to \begin{align}
\label{ch6-eq15}
&\sum_{i=0}^N\lt(\sum_{j=0}^{n-1}(i10^{-j}-\lf i10^{-j}\rf)\rt)^2\notag\\
  =&\lt(\frac{n^2}{4}-\frac{49}{108}n\rt)N-\frac{1}{2}\sum_{j=0}^{n-1}j\sum
  _{l=0}^{j-1}a_l10^l+\frac{1}{2}\sum_{j=0}^{n-1}j10^{-j}\lt(\sum_{l=0}^{j-1}
  a_l10^l\rt)^2+\cO(N).
\end{align}
The rest of the proof consists of choosing an appropriate $i_0$ (or equivalently, $y=\lala\log_{10}i_0\rara$) to obtain a sufficiently precise bound for \eqref{ch6-eq10}: Let
$i_0=
10^{n-1}-10^{\lf n/2\rf-1}+1$ if $10^{n-1}\le N\le10^n-10^{\lf n/2\rf}$, and $i_0=10^n-10^{\lf n/2\rf}$ if $10^n-10^{\lf n/2\rf}<N\le10^n-1$.  Note that $\lf N/10\rf+1\le i_0\le N-1$, and it is straightforward to verify that
        \eqb\label{ch6-eq16}
   \sum_{j=0}^{n-1}(i_010^{-j}-\lf i_010^{-j}\rf)=\frac{n}{2}+c+\cO\lt(N^{-1/2}\rt).
    \eqe
for some finite constant $c.$ Combining \eqref{ch6-eq12}, \eqref{ch6-eq15} and \eqref{ch6-eq16} yields
\begin{align*}&\sum_{i=0}^N\lt\{\sum_{j=0}^{n-1}\lt((i10^{-j}-\lf i10^{-j}\rf)-(i_010^{-j}-\lf i_010^{-j}\rf)\rt)\rt\}^2\\
=&\frac{11}{108}nN+\frac{1}{2}\sum_{j=0}^{n-1}(n-j)\lt(\sum_{l=0}^{j-1}a_l10^l
\rt)\lt(1-\sum_{l=0}^{j-1}a_l10^{l-j}\rt)+\cO(N).\end{align*}

Next, observe that \[
\frac{1}{2}\sum_{j=0}^{n-1}(n-j)\lt(\sum_{l=0}^{j-1}a_l10^l\rt)\lt(1-
\sum_{l=0}^{j-1}a_l10^{l-j}\rt)=\cO(N).
\]
which implies that \[\sum_{i=0}^N\lt\{\sum_{j=0}^{n-1}\lt((i10^{-j}-\lf i10^{-j}\rf)-(i_010^{-j}-\lf i_010^{-j}\rf)\rt)\rt\}^2=\frac{11}{108}nN+\cO(N),
\]
and hence\eqb\label{ch6-eq17}
\sum_{i=\lf N/10\rf}^N\lt\{\sum_{j=0}^{n-1}\lt((i10^{-j}-\lf i10^{-j}\rf)-(i_010^{-j}-\lf i_010^{-j}\rf)\rt)\rt\}^2
\le\frac{11\log N}{108\log10}N+\cO(N).
\eqe Let $y=\lala\log_{10}i_0\rara$. Combining \eqref{ch6-eq7}, \eqref{ch6-eq8}, \eqref{ch6-eq9} and \eqref{ch6-eq17} yields
\begin{align*}
&\int_0^1\Bigl|\bigl(F_{\nu_N\circ\iota^{-1}}(x)-F_{\eta_N\circ\iota^{-1}}(x)
\bigr)-\bigl(F_{\nu_N\circ\iota^{-1}}(y)
-F_{\eta_N\circ\iota^{-1}}(y)\bigr)\Bigr|{\rm d}x\\ &\le\frac{1}{6\log10}\sqrt{\frac{33}{\log10}}\frac{\sqrt{\log N}}{N}+\cO\lt(N^{-1}\rt);
\end{align*}and hence with Proposition \ref{ch6-prop1}, it follows at long last that \[
\limsup_{N\to\infty}\frac{N}{\sqrt{\log N}}d_{\mathbb{T}}\lt(\nu_N,\eta_N\rt)\le\displaystyle
\frac{1}{6\log10}\sqrt{\frac{33}{\log10}}.
\]
\pre
Note that Theorem \ref{ch6-thm1} describes the asymptotics of $\lt(\nu_N(\log_bn)\rt)_{N\ge1}$, in that it not only gives the rate of convergence, but also identifies the exponential distribution with specific rotation that $\lt(\nu_N\rt)$ asymptotically approaches.
\rb\label{ch6-rem1}
(i) It follows from a general result in \cite{X18} that \[
\limsup_{N\to\infty}\frac{N}{\log N}d_{\mathbb{T}}\lt(\nu_N,E_{\log b}\circ R^{-1}_{-\log_bN}\rt)<+\infty\] for every $b\in\N\setminus\{1\}$. Obviously, this is weaker than \eqref{ch6-eq4}.\\
(ii) From Zador's theorem on asymptotic quantization error in $\cP(\cT)$ \cite[Thm.1.4]{I16}, it follows that \[\liminf_{N\to\infty}Nd_{\mathbb{T}}\lt(\nu_N,E_{\log b}\circ R^{-1}_{-\log_bN}\rt)>0.\] This shows that $\lt(d_{\mathbb{T}}\lt(\nu_N,E_{\log b}\circ R^{-1}_{-\log_bN}\rt)\rt)_{N\ge1}$ cannot decay faster than\\ \noindent$\lt(N^{-1}\rt)$, and \cite[Cor.3.8]{CM95} suggests that it may be challenging to improve this lower bound.\\
(iii)
Even if the inequality \eqref{ch6-eq8} is replaced by the following H\"{o}lder inequality,
\begin{align*}
     &\sum_{i=\lf N/10\rf}^N\frac{1}{i}\lt|\sum_{j=0}^{n-1}\lt((i_010^{-j}-\lf i_010^{-j}\rf)-(i10^{-j}-\lf i10^{-j}\rf)\rt)\rt|\\
     \le&\lt(\sum_{i=\lf N/10\rf}^N\frac{1}{i^{4/3}}\rt)^{3/4}\cdot\lt(\sum_{i=\lf N/10\rf}^N\lt(\sum_{j=0}^{n-1}\lt((i_010^{-j}-\lf i_010^{-j}\rf)-(i10^{-j}-\lf i10^{-j}\rf)\rt)\rt)^4\rt)^{1/4},
    \end{align*}
 the upper bound for the rate of convergence does not improve. Indeed, a tedious computation similar to the one in the proof of Theorem \ref{ch6-thm1} yields
\eqb\label{ch6-eq19}
\limsup_{N\to\infty}\frac{N}{\sqrt{\log N}}d_{\mathbb{T}}\lt(\nu_N,E_{\log b}\circ R^{-1}_{-\log_bN}\rt)\le c,\eqe
where the constant $c$ may be smaller than $\displaystyle\frac{1}{6\log10}\sqrt{\frac{33}{\log10}}$ but still is positive.
From this, one may optimistically conjecture that for all $b>1$ (not necessarily integers), the sequence \[\lt(\frac{N}{\sqrt{\log N}}d_{\mathbb{T}}\lt(\nu_N,E_{\log b}\circ R^{-1}_{-\log_bN}\rt)\rt)_{N\ge2}\] is bounded above and below by positive constants. Especially for non-integer $b$, this is speculation only, since many of the explicit calculations and estimates leading to \eqref{ch6-eq4} do not apply directly.
\re
\subsection{Sharp rates of convergence w.r.t.\ $d_{\I}$ and $d_{\sf K}$}\label{ch6-sec2-sub2}
In this subsection, we complement the results of Subsection \ref{ch6-sec2-sub1} by characterizing the sharp rate of convergence of $\lt(\nu_N(\log_bn)\rt)_{N\ge1}$ w.r.t.\ both $d_{\I}$ and $d_{\sf K}$.

\thmb\label{ch6-thm2}
Assume $b\in\N\setminus\{1\}$. Then
\[
\lim_{N\to\infty}\frac{N}{\log N}d_{\I}\lt(\nu_N\circ\iota^{-1},E_{\log b}\circ R^{-1}_{-\log_bN}\circ\iota^{-1}\rt)=\frac{1}{2\log b}.
\]
\thme
\prb Recall that $b=10$.
As in the proof of Theorem \ref{ch6-thm1}, by formula \eqref{ch6-eq2}, it is easy to verify that for $10^{n-1}\le N\le10^n-1,$
\eqb\label{ch6-eq20}
d_{\I}\lt(\nu_N\circ\iota^{-1},E_{\log b}\circ R^{-1}_{-\log_bN}\circ\iota^{-1}\rt)
     =\frac{1}{N\log10}\sum_{i=\lf N/10\rf}^N\frac{1}{i}\sum_{j=0}^{n-1}\lt(i10^{-j}-\lf i10^{-j}\rf\rt)+\cO\lt(N^{-1}\rt).
\eqe
Like the expression for $\sum_{i=0}^N\sum_{j=0}^{n-1}(i10^{-j}-\lf i10^{-j}\rf)$ as in the proof of Theorem \ref{ch6-thm2}, it is readily checked that
\begin{align*}
  &\sum_{i=\lf N/10\rf}^N\frac{1}{i}\sum_{j=0}^{n-1}(i10^{-j}-\lf i10^{-j}\rf)\\
    =&\sum_{j=0}^{n-1}\lt(\sum_{l=\lf N/10\rf-\lf\lf N/10\rf10^{-j}\rf10^j}^{10^j-1}\frac{l10^{-j}}{\lf\lf N/10\rf10^{-j}\rf10^j+l}\rt.\\ &\lt.\quad\quad+\sum_{k=\lf\lf N/10\rf10^{-j}\rf+1}^{\lf N10^{-j}\rf-1}\sum_{l=0}^{10^j-1}\frac{l10^{-j}}{k10^j+l}+\sum_{l=0}^{N-
    \lf N10^{-j}\rf10^j}\frac{l10^{-j}}{\lf N10^{-j}\rf10^j+l}\rt),
\end{align*}
which implies that
\eqb\label{ch6-eq21}
\begin{split}
  &\sum_{i=\lf N/10\rf}^N\frac{1}{i}\sum_{j=0}^{n-1}(i10^{-j}-\lf i10^{-j}\rf)\\
  \ge&\sum_{j=0}^{n-1}\lt(\sum_{l=\lf N/10\rf-\lf\lf N/10\rf10^{-j}\rf10^j}^{10^j-1}\frac{l10^{-j}}{\lf\lf N/10\rf10^{-j}\rf10^j+10^j}\rt.\\&\lt.+\sum_{k=\lf\lf N/10\rf10^{-j}\rf+1}^{\lf N10^{-j}\rf-1}\sum_{l=0}^{10^j-1}\frac{l10^{-j}}{k10^j+10^j}+\sum_{l=0}^{N-
  \lf N10^{-j}\rf10^j}\frac{l10^{-j}}{N}\rt)\\
  =&\frac{1}{2}\sum_{j=0}^{n-1}\biggl(\frac{1-10^{-j}+\lf N/10\rf10^{-j}-\lf\lf N/10\rf10^{-j}\rf}{\lf\lf N/10\rf10^{-j}\rf+1}\cdot\lt(1-\lf N/10\rf10^{-j}+\lf\lf N/10\rf10^{-j}\rf\rt)\\ &+\frac{N10^{-j}-\lf N10^{-j}\rf}{N}\cdot\lt(N10^{-j}-\lf N10^{-j}\rf+1\rt)10^j+(1-10^{-j})\sum_{k=\lf\lf N/10\rf10^{-j}\rf+1}^{\lf N10^{-j}\rf-1}\frac{1}{k+1}\biggr).
\end{split}
\eqe
Note that
\begin{align*}
&\frac{1}{2}\sum_{j=0}^{n-1}\biggl(\frac{1-10^{-j}+\lf N/10\rf10^{-j}-\lf\lf N/10\rf10^{-j}\rf}{\lf\lf N/10\rf10^{-j}\rf+1}\cdot(1-\lf N/10\rf10^{-j}+\lf\lf N/10\rf10^{-j}\rf)\\&\quad\quad\ +\frac{(N10^{-j}-\lf N10^{-j}\rf)}{N}\cdot(N10^{-j}-\lf N10^{-j}\rf+1)10^j\biggr)=\cO(1).\end{align*}
Moreover, it is tedious but straightforward to confirm that \[\frac{1}{2}\sum_{j=0}^{n-1}(1-10^{-j})\sum_{k=\lf\lf N/10\rf10^{-j}\rf+1}^{\lf N10^{-j}\rf-1}\frac{1}{k+1}
=\frac{1}{2}\log N+\cO(1),
\]
with which \eqref{ch6-eq21} takes the form \eqb\label{ch6-eq22}
 \sum_{i=\lf N/10\rf}^{N}\frac{1}{i}\sum_{j=0}^{n-1}(i10^{-j}-\lf i10^{-j}\rf)\ge \frac{1}{2}\log N+\cO(1).\eqe
Analogously, one can also show that \eqref{ch6-eq22} holds with $\ge$ replaced by $\le$, and hence\[\sum_{i=\lf N/10\rf}^N\frac{1}{i}\sum_{j=0}^{n-1}(i10^{-j}-\lf i10^{-j}\rf)=\frac{1}{2}\log N+\cO(1).\] The conclusion now follows from \eqref{ch6-eq20}.
\pre
\noindent
The following corollary is immediately obtained from Theorem \ref{ch6-thm2}, together with \cite[Thm.5]{OS16} and the fact that $d_{\I}\le d_{\sf K}$.
\cob\label{ch6-cor1}
Assume $b\in\N\setminus\{1\}$. Then
\begin{align*}
0<&\liminf_{N\to\infty}\frac{N}{\log N}d_{\sf K}\lt(\nu_N\circ\iota^{-1},E_{\log b}\circ R^{-1}_{-\log_bN}\circ\iota^{-1}\rt)\\
\le&\limsup_{N\to\infty}\frac{N}{\log N}d_{\sf K}\lt(\nu_N\circ\iota^{-1},E_{\log b}\circ R^{-1}_{-\log_bN}\circ\iota^{-1}\rt)<+\infty.
\end{align*}
\coe
Comparing Theorems \ref{ch6-thm1} and \ref{ch6-thm2}, as well as Corollary \ref{ch6-cor1}, notice how\\ \noindent $\lt(d_{\cT}\lt(\nu_N,E_{\log b}\circ R^{-1}_{-\log_bN}\rt)\rt)_{N\ge1}$ decays somewhat faster than both\\ \noindent $\lt(d_{\I}\lt(\nu_N\circ\iota^{-1},E_{\log b}\circ R^{-1}_{-\log_bN}\circ\iota^{-1}\rt)\rt)_{N\ge1}$ and $\lt(d_{\sf K}\lt(\nu_N\circ\iota^{-1},E_{\log b}\circ R^{-1}_{-\log_bN}\circ\iota^{-1}\rt)\rt)_{N\ge1}$. Moreover, the ratio \[ \frac{d_{\I}\lt(\nu_N\circ\iota^{-1},E_{\log b}\circ R^{-1}_{-\log_bN}\circ\iota^{-1}\rt)}{d_{\sf K}\lt(\nu_N\circ\iota^{-1},E_{\log b}\circ R^{-1}_{-\log_bN}\circ\iota^{-1}\rt)}\] is bounded above and below by positive constants. This is remarkable since
$$\inf\nolimits_{\mu\neq\nu,\ \mu,\nu\in\cP(\cT)}\frac{d_{\I}\lt(\mu\circ\iota^{-1},\nu\circ\iota^{-1}\rt)}
{d_{\sf K}\lt(\mu\circ\iota^{-1},\nu\circ\iota^{-1}\rt)}=0.$$

\subsection*{Acknowledgements}

Deepest thanks to the author's PhD thesis advisor Arno Berger at the University of Alberta for proposing this project, many helpful discussions, as well as the effort in improving the presentation of this manuscript as one chapter in the author's thesis. The author is also indebted to one anonymous referee for his/her serious work including pointing out a mistake in the proof of Theorem~\ref{ch6-thm1} as well as other valuable comments which also help improve the presentation of the current manuscript. This research is supported in part by a Pacific Institute for the Mathematical Sciences (PIMS) Graduate Scholarship and a Josephine Mitchell Graduate Scholarship at the University of Alberta.

\end{document}